\documentclass[12pt]{article}
\usepackage{color}
\usepackage{amssymb, amsmath}
\usepackage{epsfig}
\usepackage{graphicx}

\newcommand{\bth}{\begin{theorem}}    \newcommand{\ethe}{\end{theorem}}
\newcommand{\bre}{\begin{remark}\em }  \newcommand{\ere}{\end{remark}}
\newcommand{\ble}{\begin{lemma}}       \newcommand{\ele}{\end{lemma}}
\newcommand{\bde}{\begin{definition}}   \newcommand{\ede}{\end{definition}}
\newcommand{\bco}{\begin{corollary}}     \newcommand{\eco}{\end{corollary}}
\newcommand{\bpr}{\begin{proposition}}    \newcommand{\epr}{\end{proposition}}
\newcommand{\bexer}{\begin{exercise}}       \newcommand{\eexer}{\end{exercise}}
\newcommand{\bexam}{\begin{example}}        \newcommand{\eexam}{\end{example}}
\newcommand{\bfi}{\begin{fig}}                  \newcommand{\efi}{\end{fig}}
\newcommand{\btab}{\begin{tab}}         \newcommand{\etab}{\end{tab}}
\newcommand{\bpf}{\begin{proof}}            \newcommand{\epf}{\end{proof}}
\newcommand{\bce}{\begin{center}}   \newcommand{\ece}{\end{center}}

\newcommand{\D}{\Delta}

\newcommand{\bfE}{{\bf E}}
\newcommand{\bbf}{{\mathcal F}}

\newcommand{\bfP}{{\bf P}}

\newcommand{\barr}{\begin{array}}   \newcommand{\earr}{\end{array}}
\newcommand{\beao}{\begin{eqnarray*}}   \newcommand{\eeao}{\end{eqnarray*}\noindent}
\newcommand{\beam}{\begin{eqnarray}}    \newcommand{\eeam}{\end{eqnarray}\noindent}
\newcommand{\beqq}{\begin{equation}}    \newcommand{\eeqq}{\end{equation}\noindent}

\newtheorem{theorem}{Theorem}
 \newtheorem{corollary}{Corollary}
 \newtheorem{proposition}{Proposition}
\newtheorem{lemma}{Lemma}%

\newtheorem{remark}{Remark}

\newtheorem{definition}{Definition}
\newtheorem{example}{Example}
\newcommand{\proof}{\noindent\textbf{Proof.~}}

\newcommand{\qed}{\space\hfill\hspace*{\fill} $\vbox{\hrule\hbox{\vrule
height1.3ex\hskip1.3ex\vrule}\hrule}$\hss\vskip\topsep\relax}

\begin{document}

\title{Construction of positivity preserving numerical schemes for multidimensional stochastic differential equations}

\author{Nikolaos Halidias \\
{\small\textsl{Department of Mathematics }}\\
{\small\textsl{University of the Aegean }}\\
{\small\textsl{Karlovassi  83200  Samos, Greece} }\\
{\small\textsl{email: nikoshalidias@hotmail.com}}}

\maketitle

\begin{abstract} In this note we work on the construction of positive
preserving numerical schemes for systems of stochastic
differential equations. We use the semi discrete idea that we have
proposed before  proposing now a numerical scheme that preserves
positivity on multidimensional stochastic differential equations
converging strongly in the mean square sense to the true solution.
\end{abstract}



\section{Introduction}
Throughout, let $T>0$ and $(\Omega, \bbf, \{\bbf_t\}_{0\leq t\leq
T}, \bfP)$ be a complete probability space, meaning that the
filtration $ \{\bbf_t\}_{0\leq t\leq T} $ satisfies the usual
conditions, i.e. is right continuous and $\bbf_0$ includes all
$\bfP-$null sets, and let  an $m$-dimensional Wiener process
$W(t)$ defined on this space.

Let the following multidimensional stochastic differential
equation,
\begin{eqnarray}
x(t) = x(0) + \int_0^t a(x(s)) ds  +  \sum_{j=1}^m \int_0^t
b_j(x(s))dW_j(s),
\end{eqnarray}
where $a(\cdot) : \mathbb{R}^d \to \mathbb{R}^d$, $b_j(\cdot) :
\mathbb{R}^{d} \to  \mathbb{R}^{d}, j=1,..,m$, $x(0) = x$ and
$x:\Omega \to \mathbb{R}^d$ is ${\cal F}_0$ measurable.  That is,
$a(x) = (a_1(x),...,a_d(x))$, $b_j(x) = (b_{1},...,b_{d})$ where
$j = 1,...,m$.

\textbf{Assumption A}  Suppose that there exists some functions
$f_i(x,y) : \mathbb{R}^{2d} \to \mathbb{R}$ and $g_{ij}(x,y) :
\mathbb{R}^{2d} \to \mathbb{R}$ such that $f_i(x,x) = a_i(x)$ and
$g_{ij}(x,x) = b_{ij}(x)$ with $x = (x_1,...,x_d)$, $y =
(y_1,...,y_d)$, $i = 1,...,d$ and $j=1,...,m$. Let $\mathbb{E}
\|x\|_2^{p} < A$  for some $p
> 2$. Suppose further that $f,g$ satisfy the following condition
$$
|f_i(x_1,y_1) - f_i(x_2,y_2)| + |g_{ij}(x_1,y_1) -
g_{ij}(x_2,y_2)| \leq C_R \left(\|x_1-x_2\|_2 +
\|y_1-y_2\|_2\right),
$$
for  any $R>0$, $i=1,...,d$, $j=1,...,m$ and $x_1,x_2,y_1,y_2$
such that $\|x_1\|_2\vee\|x_2\|_2\vee\|y_1\|_2\vee\|y_2\|_2\leq
R,$ where the constant $C_R$ depends on $R$ and $x\vee y$ denotes
the maximum of $x, y.$ Here  $$ \| x \|_p = \sqrt[p]{\sum_{i=1}^d
x_i^p}
$$
where $x = (x_1,...,x_d)$.

 Let the equidistant partition
$0=t_0<t_1<...<t_N=T$ and $\Delta=T/N.$ We propose the following
numerical scheme,
\begin{eqnarray}
y_{i}(t) = x_i + \int_0^t f_i(y(s),y(\hat{s})) ds +
 \sum_{j=1}^m \int_0^t
 g_{ij}(y(s),y(\hat{s})) dW_j(s).
\end{eqnarray}
for $i=1,...,d$. Note that $y(s) = (y_1(s),...,y_d(s))$ and
$y(\hat{s}) = (y_1(\hat{s}),...,y_d(\hat{s}))$, with
$y_{i}(\hat{s}) = y_{i}(t_k)$ when $s \in [t_k,t_{k+1}]$ and
$y_i(0) = x_i$. This is again, in general, a system of stochastic
differential equations and we suppose that has a unique strong
solution. However, in practice, we will choose $f,g$ in  a way
that the resulting system has less than $d$ dependent equations
or/and having known explicit solution. An interesting choice that
we will see in our example is that the resulting system is not in
fact a system of SDEs but $d$
 independent equations. Each of these equations is linear with
known explicit solution and we solve it independently of the
others.

Below we state a convergence, in the mean square sense, theorem of
$y_t$ to the true solution   as $\Delta \downarrow 0$. The proof
of this theorem is exactly the same as in \cite{halidias:2013}
changing the absolute values by the Euclidean norm in
$\mathbb{R}^d$.

\begin{theorem} Suppose Assumption A holds and (2) has a unique
strong solution.  Let also
$$
\bfE(\sup_{0\leq t\leq T}\|x_t\|^p_2) \vee \bfE(\sup_{0\leq t\leq
T}\|y_t\|_2^p)<A,
$$
for some $p>2$ and $A>0.$ Then the semi-discrete numerical scheme
(2) converges to the true solution of (1) in the mean square
sense, that is \beqq \label{eq600}
\lim_{\D\rightarrow0}\bfE\sup_{0\leq t\leq T}\|y_t-x_t\|_2^2=0.
\eeqq
\end{theorem}

\proof We use the same arguments as in \cite{halidias:2013},
\cite{Mao}. Set $ \rho_R = \inf \{ t \in [0,T]:
 \|x(t)\|_2 \geq R \}$ and $ \tau_R = \inf \{ t \in [0,T] : \|y(t)\|_2 \geq R
\}$. Let $\theta_R = \min \{ \tau_R, \rho_R \}$. Using exactly the
same arguments as in \cite{Mao} we obtain,
\begin{eqnarray*}
\mathbb{E}\left(\sup_{0 \leq t \leq T } \|y(t) -x(t)\|_2^2\right)
\leq \mathbb{E}\left( \sup_{ 0 \leq t \leq T} \left\|y({t
\wedge\theta_R}) - x({t \wedge\theta_R})\right\|_2^2\right) +
\frac{2^{p+1} \delta A}{p} + \frac{(p-2)2A}{p
\delta^{\frac{2}{p-2}}R^p}.
\end{eqnarray*}

First, let us estimate the  quantity $\mathbb{E} \|y({t\wedge
\theta}) - y({\widehat{t\wedge \theta}})\|_2^2 $, beginning with,
\begin{eqnarray*}
|y_i({t\wedge \theta}) - y_i({\widehat{t\wedge \theta_R}})|^2 =
\left| \int_{{\widehat{t\wedge \theta_R}}}^{t\wedge \theta_R}
f_i(y_s,y_{\hat{s}})ds + \sum_{j=1}^m \int_{{\widehat{t\wedge
\theta_R}}}^{t\wedge \theta_R} g_{ij}( y_s,y_{\hat{s}})dW_j(s)
\right|^2 \leq
\\ C \left( \left(\int_{{\widehat{t\wedge \theta_R}}}^{t\wedge \theta_R} f_i(y_s,y_{\hat{s}})ds \right)^2 +
\sum_{j=1}^m \left|\int_{{\widehat{t\wedge \theta_R}}}^{t\wedge
\theta_R} g_{ij}(y_s,y_{\hat{s}})dW_j(s) \right|^2 \right).
\end{eqnarray*}
Taking expectations, using Ito's isometry and  the fact that
$|f_i(y_s,y_{\hat{s}})|, | g_{ij}(y_s,y_{\hat{s}})| \leq C_R$ we
have that,
\begin{eqnarray*}
\mathbb{E} |y_i({t\wedge \theta}) - y_i({\widehat{t\wedge
\theta}})|^2 \leq C_R \Delta,
\end{eqnarray*}
and from this it follows that
\begin{eqnarray*}
\mathbb{E} \|y({t\wedge \theta}) - y({\widehat{t\wedge
\theta}})\|_2^2 \leq C_R \Delta,
\end{eqnarray*}

Next we work on the quantity $\mathbb{E}\left( \sup_{ 0 \leq t
\leq T} \left\|y({t \wedge\theta_R}) - x({t
\wedge\theta_R})\right\|_2^2\right)$ to get
\begin{eqnarray*}
& & |x_i({t\wedge \theta_R}) - y_i({t\wedge \theta_R})|^2 = \\
 & & \left| \int_0^{t\wedge \theta_R}
(f_i(x_s,x_s)-f_i(y_s,y_{\hat{s}}))  ds + \sum_{j=1}^m
\int_0^{t\wedge \theta_R}  (g_{ij}(x_s,x_s)-g_{ij}(y_s,y_{\hat{s}}))dW_j(s) \right|^2 \\
& & \leq C \left( \int_0^{t\wedge \theta_R}
|f_i(x_s,x_s)-f_i(y_s,y_{\hat{s}})|^2   ds +  \sum_{j=1}^m
\left|\int_0^{t\wedge
\theta_R} (g_{ij}(x_s,x_s)-g_{ij}(y_s,y_{\hat{s}}))dW_j(s) \right|^2 \right) \\
& & \leq C_R \int_0^{t\wedge \theta_R} (\|x_s-y_s\|_2^2 +
\|x_s-y_{\hat{s}}\|_2^2)ds + C \sum_{j=1}^m\left|\int_0^{t\wedge
\theta_R}
(g_{ij}(x_s,x_s)-g_{ij}(y_s,y_{\hat{s}}))dW_j(s)\right|^2
\end{eqnarray*}
 We can write
now,
\begin{eqnarray*}
  \sup_{0 \leq t \leq s} |x_i({t\wedge \theta_R}) - y_i({t\wedge
\theta_R})|^2 \end{eqnarray*}
\begin{eqnarray*}
   \leq & & C_R \int_0^s \left(\|x({r\wedge
\theta_R})-y({r\wedge \theta_R})\|_2^2 + \|x({r\wedge \theta_R}) -
y({\widehat{r\wedge \theta_R}})|_2^2 \right)dr  \\ & &  + C
\sum_{j=1}^m \sup_{0 \leq t \leq s} \left|\int_0^{t\wedge
\theta_R}
(g_{ij}(x_s,x_s)-g_{ij}(y_s,y_{\hat{s}}))dW_j(s)\right|^2.
\end{eqnarray*}
Taking expectations on both sides, using  Doob's martingale
inequality for the second term at the right hand side and
Assumption A for $g(\cdot,\cdot)$ we arrive at,
\begin{eqnarray*}
\mathbb{E} (\sup_{0 \leq t \leq s} |x_i({t\wedge \theta_R}) -
y_i({t\wedge \theta_R})|^2) & \leq & C_R \mathbb{E} \int_0^s
\big(\|x({r\wedge \theta_R})-y({r\wedge \theta_R})\|_2^2  \\ & & +
\|x({r\wedge \theta_R}) - y({\widehat{r\wedge \theta_R}})\|_2^2
\big)dr
\\ & \leq & C_R \int_0^s \mathbb{E} \sup_{ 0 \leq l \leq r} \|x({l\wedge
\theta_R})-y({l\wedge \theta_R})\|_2^2 dr \\ & & + C_R \int_0^s
\mathbb{E} \|y({r\wedge \theta_R}) - y({\widehat{r\wedge
\theta_R}})\|_2^2 dr.
\end{eqnarray*}

Therefore, we have that
\begin{eqnarray*}
\mathbb{E}\left( \sup_{ 0 \leq t \leq T} \left\|y({t
\wedge\theta_R}) - x({t \wedge\theta_R})\right\|_2^2\right) & \leq
& \mathbb{E} \sum_{i=1}^d (\sup_{0 \leq t \leq s} |x_i({t\wedge
\theta_R}) - y_i({t\wedge \theta_R})|^2) \\ & \leq  & C_R \int_0^s
\mathbb{E} \sup_{ 0 \leq l \leq r} \|x({l\wedge
\theta_R})-y({l\wedge \theta_R})\|_2^2 dr + C_R \Delta.
\end{eqnarray*}
Using Gronwall's inequality and continuing as in
\cite{halidias:2013} we arrive at the desired result.\qed

\section{Example}

Consider the following multidimensional SDE.

\begin{eqnarray*}
x(t) = x + \int_0^t \left(x(s) - \| x(s) \|_2^2 x(s) \right)ds +
\int_0^t x(s) dW(s).
\end{eqnarray*}

{\bf Assumption B} Assume that $\mathbb{E} | \ln x_i | +
\mathbb{E} x_i^{2p} <C$ for all $i=1,...,d$ for some $p > 2$ and
$C > 0$.

Consider the following system of SDEs,
\begin{eqnarray*}
y_{i}(t) = x_{i} + \int_0^t \left(y_{i}(s) - \| y({\hat{s}})
\|_2^2 y_i(s)\right) ds + \int_0^t y_{i}(s) dW(s).
\end{eqnarray*}
Each of the above equation has only one unknown function and is
linear with known explicit solution and thus preserves positivity.

\begin{theorem} Under Assumption B  we have
\begin{eqnarray*}
\mathbb{E}(\sup_{0 \leq t \leq T} \| x_t \|_2^p) < A, \quad
\mathbb{E}(\sup_{0 \leq t \leq T} \| y_t \|_2^p) < A
\end{eqnarray*}
for some $A > 0$.
\end{theorem}

\proof For the moment bound of the true solution we use Lemma 3.2
of  \cite{Mao} since the drift coefficient satisfies the
monotonicity condition and the diffusion term the linear growth
condition.

We will prove now that the approximate solution has bounded
moments. Set the stopping time $\theta_R = \inf \{ t \in [0,T]:
y_i(t) > R \}$. Using Ito's formula on $y_i^q(t \wedge \theta_R)$
(with $q = 2p$)
 we obtain
\begin{eqnarray*}
y_i^q(t \wedge \theta_R) = x_i^q + \int_0^{t \wedge \theta_R}
y_i^q(s) \left(q - \| y(\hat{s}) \|_2^2 +\frac{q(q-1)}{2} \right)
ds + \int_0^{t \wedge \theta_R} q y_i^q(s) dW(s).
\end{eqnarray*}
Taking expectations we arrive at the following inequality,
\begin{eqnarray*}
\mathbb{E} y_i^q(t \wedge \theta_R) \leq \mathbb{E} x_i^q
+\left(q+\frac{q(q-1)}{2}\right) \int_0^{t \wedge \theta_R}
\mathbb{E} y_i^q(s) ds.
\end{eqnarray*}
Using now Gronwall's inequality we get that $\mathbb{E} y_i^q(t
\wedge \theta_R) < A$ and using Fatou's lemma we arrive at the
bound $\mathbb{E} y_i^q(t) < A$ for $i = 1,...,d$. Therefore, we
have that $\mathbb{E}( \| y_t \|_q^q) < A$.

Using  again Ito's formula on $y_i^p(t)$  we obtain,
\begin{eqnarray*}
y_i^p(t) \leq x_i^p +\left(p + \frac{p(p-1)}{p}\right) \int_0^t
y_i^p(s)ds + \int_0^t p y_i^p(s) dW(s).
\end{eqnarray*}
Taking  the supremum over $[0,T]$ we have,
\begin{eqnarray*}
\sup_{ 0 \leq t \leq T} y_i^p(t) \leq C \left( x_i^p +\int_0^T
y_i^p(s) ds + \sup_{0 \leq t \leq T}  \int_0^t y_i^p(s) dW(s)
\right).
\end{eqnarray*}
Taking expectations and using Doob's inequality on the diffusion
term we arrive at
\begin{eqnarray*}
\mathbb{E} \left(\sup_{ 0 \leq t \leq T}  y_i^p (t) \right) \leq
A.
\end{eqnarray*}
But note that
\begin{eqnarray*}
\mathbb{E} \left(\sup_{ 0 \leq t \leq T} \| y(t) \|^p_p (t)
\right) \leq \mathbb{E} \left(\sum_{i=1}^d \sup_{ 0 \leq t \leq T}
y_i^p (t) \right) < A.
\end{eqnarray*}
Using  the equivalence of the norms in $\mathbb{R}^d$ we obtain
the desired result.
 \qed

\begin{theorem}
For the true solution we have $x_i(t) > 0$ a.s. for $i=1,...,d$
\end{theorem}

\proof Set the stopping time $\theta_R = \inf \{ t \in [0,T] :
x_i(t) < \frac{1}{R} \}$.
 Using Ito's formula on $\ln x_i(t \wedge \theta_R)$ we obtain
\begin{eqnarray*}
\ln x_i(t \wedge \theta_R) = \ln x_i  & + & \int_0^{t \wedge
\theta_R} \left( \frac{1}{x_i(s)} \left(x_i(s) - \| x(s) \|_2^2
x_i(s) \right) - \frac{1}{2} \frac{1}{x_i^2(s)} x_i^2(s) \right) ds \\
& + & \int_0^{t \wedge \theta_R} \frac{1}{x_i(s)} x_i(s) dW(s).
\end{eqnarray*}
Taking absolute values and then expectations, using Jensen
inequality and then Ito's isometry on the diffusion term we arrive
at
\begin{eqnarray*}
\mathbb{E}| \ln x_i(t \wedge \theta_R) | \leq \mathbb{E} | \ln x_i
| + \int_0^{t \wedge \theta_R} \left(1 + \mathbb{E} \| x(s)
\|_2^2\right) ds +\sqrt{T} \leq C.
\end{eqnarray*}
We have used also the moment bound for the true solution  (see
Theorem 2).
 Therefore
 \begin{eqnarray*}
 \ln R \;\;P\left(\{\theta_R \leq R\}\right) < C,
 \end{eqnarray*}
 and thus $P\left(\{\theta_R \leq R\}\right) \to 0$ as $R \to \infty$.
 But
 \begin{eqnarray*}
 P\left(\{x_i(t) \leq 0\}\right) =  P\left(\bigcap_{R=1}^{\infty} \{ x_i(t) < \frac{1}{R} \}\right) =
\lim_{R \to \infty} P\left(\{ x_i(t) < \frac{1}{R} \}\right) \leq
\lim_{R \to \infty} P\left(\{\theta_R < t\}\right)=0.
\end{eqnarray*}
 \qed

In this example we have $a(x) = x - \|x\|_2^2 x$ and $b(x) = x$
with $x = (x_1,...,x_d)$, that is $m=1$. We choose $f_i(x,y) = x_i
- \| y \|_2^2 x_i$ and $g_i(x,y) = x_i$ where $x = (x_1,...,x_d)$
and $y = (y_1,...,y_d)$. It is easy to see that $f_i(x,x) = a(x)$
and $g_i(x,x) = b(x)$. Moreover, $f_i,g_i$ satisfies Assumption A.
Therefore, our proposed numerical scheme preserves positivity and
converges strongly in the mean square sense to the true solution.

{\bf Conclusion} Concerning super linear SDEs it is well known
that the usual Euler scheme diverges and the tamed-Euler scheme
\cite{Jentzen} does not preserve positivity. For scalar SDEs there
are some numerical schemes that preserves positivity (see for
example
\cite{halidias:2013},\cite{higham_et_al:2012},\cite{Liu},\cite{neuenkirch_szpruch:2012}
and the references therein) but for the multidimensional case it
is not clear  how can be extended. So we extend here our semi
discrete method (\cite{halidias:2013}), that we have proposed
before for scalar SDEs, to the multidimensional case. There is
also the possibility in some multidimensional SDEs to combine the
semi discrete method with another method designed for scalar SDEs
in order to construct positivity preserving numerical schemes. Let
us note that we can use the semi discrete method also in the case
when the diffusion term is super linear.
 Our goal in the future is to  apply the semi discrete method
to more complicated systems of SDEs.

\end{document}